\documentclass[12pt]{article}
\usepackage{amsmath,amsthm,amsfonts,amssymb,amscd}
\usepackage{graphicx}
 \normalsize
\usepackage[cp1251]{inputenc}  
\usepackage[T2B]{fontenc}      
\usepackage{color}

\newtheorem{theorem}{Theorem}[section]
\newtheorem{lemma}[theorem]{Lemma}
\newtheorem{prop}[theorem]{Proposition}

\newtheorem{rem}[theorem]{Remark}

\title{\bf Spectral equivalence of Gaussian \\
random functions: operator approach }
\author{
A.I. Nazarov\footnote{ St.Petersburg Department of Steklov Mathematical Institute of Russian Academy of Science, Fontanka 27, St.Petersburg, 191023, Russia, and St.Petersburg State University, Universitetskaya emb. 7-9, St.Petersburg, 199034, Russia; E-mail: nazarov@pdmi.ras.ru}\setcounter{footnote}{6}
\ and
\fbox{Ya.Yu. Nikitin}\footnote{St.Petersburg State University, Universitetskaya emb. 7-9, St.Petersburg, 199034, Russia, and National Research University -- Higher School of Economics, Soyuza Pechatnikov 16, St.Petersburg, 190008, Russia; E-mail: y.nikitin@spbu.ru
}
}
\date{}

\begin{document}

\maketitle

\begin{abstract}
We introduce a new approach to the spectral equivalence of Gaussian processes and fields, based on the methods of operator theory in Hilbert space. Besides several new results including identities in law of quadratic norms for integrated and multiply integrated Gaussian random functions we give an application to goodness-of-fit testing.
\end{abstract}


\medskip

\section{Introduction}

The distribution of quadratic functionals for Gaussian random functions is an interesting and intensively developing topic in connection with the demands of asymptotic problems of empirical processes (see \cite{Nik}, \cite{VNB} and references therein) and the theory of small deviations in $L_2$ (see, e.g., \cite{Lif99}, \cite{LS}). We will be interested in different Gaussian random functions with equally distributed quadratic norms.

By virtue of the Karhunen--Lo\`eve expansion, see \cite{AG}, \cite[Ch. X]{Lo}, such norms can be represented as infinite quadratic Gaussian forms, whose coefficients are the eigenvalues of the corresponding covariance operators. Therefore, to prove the identity  in law for $L_2$-norms of two Gaussian random functions, it suffices to verify the coincidence of spectra for their covariance operators (excluding zero). In this case we prefer to call two random Gaussian random functions ${\bf X}$ and ${\bf Y}$ {\it spectrally equivalent} and write ${\bf X \sim Y}$.

It should be noted that the equality of $L_2$-norms does not entail the equality of Gaussian processes or fields in law. Consider the typical example of this kind. Let $W$ be the standard Brownian motion  on $[0,1]$, and denote by $B$ the standard Brownian bridge on the same interval.

The following spectral equivalence is well known:
\begin{equation}
\label{Wa}
  \qquad W(t)- \int\limits_0^1 W (s)\,ds \sim \, B(t).
\end{equation}

For instance, Donati-Martin and Yor \cite{Dona} proved (\ref{Wa}) using the Fubini--Wiener technique, while in \cite{BNO} this equivalence was proved by the direct calculation of spectra. On the other hand, the Gaussian processes in (\ref{Wa}) have different covariances:
\begin{align*}
&\mathbb{E}\Big (W(s)- \int\limits_0^1 W (u)\,du\Big)\Big( W(t)- \int\limits_0^1 W (u)\,du\Big)
= \min(s,t) - \frac {2s-s^2}2 - \frac {2t-t^2}2 + \frac 13;\\
&\mathbb{E}B(s) B(t) = \min(s,t) - st.
\end{align*}

Peccati and Yor \cite{PecY}, Deheuvels, Peccati and Yor \cite{DPY} and Deheuvels \cite{Deh} extended spectral equivalence (\ref{Wa}) to the Gaussian fields on the unit square. Let ${\bf W}$ and ${\bf B}$ denote, respectively, the classical Brownian sheet and the bivariate Brownian bridge or pinned Brownian sheet, see \cite{DPY}. The authors of \cite{PecY} and \cite{DPY} obtained several spectral equivalences generalizing (\ref{Wa}). The simplest of them has the form:
\begin{equation}
\label{bivar}
\qquad {\bf W}(t_1,t_2)- \int\limits_0^1 \int\limits_0^1 {\bf W} (s_1,s_2) ds_1 ds_2 \sim  {\bf B}(t_1,t_2).
\end{equation}
Deheuvels, Peccati and Yor stated these results in dimension  $d= 2$, indicating that similar equivalences can be written out in the case  $d > 2$ \ ``at the price of minor additional technicalities.'' However, these formulations never appeared.

Our aim is to prove far reaching and sometimes unexpected generalizations of these and similar spectral equivalences in a very short and compact way which differs from that of \cite{DPY} and \cite{PecY}. The originality of our approach is due to the use of the methods of operator theory in Hilbert space. The spectral equivalence of various $d$-parametric Brownian functions takes on a simple and uniform perspective. We argue that the operator language is the most correct and convenient for writing identities in law  for quadratic norms of Gaussian random functions.

Moreover, we demonstrate that similar relations are also valid for {\it integrated} and {\it multiply integrated} Gaussian fields and processes, something that has been done for the first time. We expect  future applications to nonparametric statistics (in particular, to goodness-of-fit testing and testing of independence), and to the theory of Brownian functionals.

\medskip

The structure of the paper is as follows. First, we introduce the necessary facts from functional analysis in Hilbert space. Next, we represent the operations of centering, integration and bridge construction for Gaussian random functions in operator terms and prove several theorems about spectral equivalence of various processes and fields. The last section deals with the spectral properties of the kernel in the case of multivariate $\omega^2$-statistic.

\medskip

 The Gaussian processes and operators acting in the space of functions of one variable are denoted by capital letters, the multivariate Gaussian fields and the corresponding operators by bold letters. If we need to mention that, say, an operator $T$ acts on functions of variable $x_k$, we write $T_k$.

\section{Some functional analytic preliminaries}

The following statement is well known, see, e.g. \cite[Section 3.10]{BS10}.

 \begin{prop}\label{Prop1} 
 Let $A$ and $B$ be compact operators in the Hilbert space $H$. Then the non-zero eigenvalues of the operators $AB$ and $BA$ coincide (with the multiplicities).
\end{prop}
 
Recall that if $A$ and $\widetilde A$ are operators in the Hilbert spaces $H$ and $\widetilde H$ respectively, we can define their tensor product $A\otimes\widetilde A$ in the Hilbert space $H\otimes\widetilde H$. If $A$ and $\widetilde A$ are integral operators with kernels ${\cal A}(x,y)$ and $\widetilde {\cal A}(s,t)$ then $A\otimes\widetilde A$ is also an integral operator with kernel $\mathbb{A}((x,s),(y,t))={\cal A}(x,y)\widetilde {\cal A}(s,t)$.

As the eigenvalues of the tensor product $A\otimes\widetilde A$ are the products of the eigenvalues $\lambda_i(A)$ and $\lambda_j(\widetilde A)$, the following statement is obvious.

\begin{prop}\label{Prop2} 
Let $A_k$ and $B_k$ be compact operators in the Hilbert spaces $H_k$, $k=1,\dots,d$. Assume that non-zero eigenvalues of the operators $A_k$ and $B_k$ coincide (with the multiplicities) for every $k$. Then the non-zero eigenvalues of the tensor products
$$
A=\underset{k=1}{\overset{d}\otimes} A_k\qquad \text{and} \qquad B=\underset{k=1}{\overset{d}\otimes} B_k
$$
coincide (with the multiplicities).
\end{prop}

We define some operators in $L_2([0,1])$:  operators of integration from the left and from the right
$$
(Tu)(x)=\int\limits_0^{x}u(t)\,dt, \qquad (T^*u)(x)=\int\limits_{x}^1 u(t)\,dt,
$$
the orthogonal projector onto the subspace of constants, the symmetry and the multiplication operators
$$
(Pu)(x)=\int\limits_0^1 u(t)\,dt, \qquad (Su)(x)=u(1-x), \qquad (M_fu)(x)=f(x)u(x).
$$
Notice that $T^*=STS$. Also it is obvious that $PS=SP=P$ 
and $SM_fS=M_{Sf}$.

Also, we define multidimensional operators which are tensor products of the corresponding one-dimensional operators:
$$
{\bf T}=\underset{k=1}{\overset{d}\otimes} T_k, \qquad
{\bf T}^*=\underset{k=1}{\overset{d}\otimes} T_k^*, \qquad
{\bf P}=\underset{k=1}{\overset{d}\otimes} P_k, \qquad
{\bf S}=\underset{k=1}{\overset{d}\otimes} S_k. 
$$

Direct calculation shows that the covariance operator $K_W$ of the Wiener process $W(t)$ and the covariance operator $K_{W_1}$ of the inverted Wiener process $W_1(t)\equiv W(1-t)$ are given by
$$
K_W=TT^* \qquad \mbox{and} \qquad K_{W_1}=T^*T,
$$
respectively.

Furthermore, the covariance operator $K_B$ of the Brownian bridge allows for the following representation:
\begin{equation}
\label{bridge}
K_B=T(I-P)T^*=T^*(I-P)T
\end{equation}
(here $I$ stands for the identity operator).

As a corollary, we obtain representations of the covariance operators of $d$-variate Brownian sheet ${\bf W}({\bf x})$, of $d$-variate inverted Brownian sheet ${\bf W}_1({\bf x})={\bf W}({\bf 1}-{\bf x})$, and of $d$-variate Brownian pillow ${\bf B}_*({\bf x})$ (each of them is the tensor product of the corresponding Gaussian processes, see \cite{KNN}):
$$
{\bf K}_{\bf W}={\bf T}{\bf T}^*, \quad {\bf K}_{{\bf W}_1}={\bf T}^*{\bf T},\quad {\bf K}_{{\bf B}_*}=
{\bf T}\cdot\underset{k=1}{\overset{d}\otimes}(I-P_k)\cdot{\bf T}^*={\bf T}^*\cdot\underset{k=1}{\overset{d}\otimes}(I-P_k)\cdot{\bf T}.
$$

The following statements can be also easily verified by direct calculation.

\begin{lemma}\label{Le1} 
The covariance operator of $d$-dimensional pinned Brownian sheet on $[0,1]^d$
$$
{\bf B}({\bf x})={\bf W}({\bf x})-{\bf W}({\bf 1})\prod_{k=1}^d x_k
$$
has the following the representation:
$$
{\bf K}_{\bf B}={\bf T}({\bf I}-{\bf P}){\bf T}^*.
$$
\end{lemma}

\begin{rem}\label{Rem} 
 Notice that for $d>1$ the function ${\bf B}$ is not 
symmetric, and thus we have no multivariate counterpart 
of the last identity in (\ref{bridge}).
\end{rem}

\begin{lemma}\label{Le2} 
Let the Gaussian field $\bf X$ on $[0,1]^d$ have the covariance operator ${\bf K}_{\bf X}$. Then

1. The covariance operator of the inverted field
${\bf X}({\bf 1}-{\bf x})$ has the following representation:
$$
{\bf K}_{{\bf X}({\bf 1}-{\bf x})}={\bf S}{\bf K}_{\bf X}{\bf S}.
$$

2. The covariance operator of the centered field
$$
\overline{\bf X\vphantom{^1}}({\bf x})={\bf X}({\bf x})-\int\limits_{[0,1]^d}{\bf X}({\bf y})\,d{\bf y}
$$
has the following representation:
$$
{\bf K}_{\overline{\bf X\vphantom{^1}}}=({\bf I}-{\bf P}){\bf K}_{\bf X}({\bf I}-{\bf P}).
$$

3. The covariance operators of the (left- and right-) integrated fields
$$
{\bf X}^{[0]}({\bf x})=\int\limits_0^{x_1}\dots\int\limits_0^{x_d}{\bf X}({\bf  y})\,dy_1\dots dy_d \quad \text{and}\quad
{\bf X}^{[1]}({\bf x})=\int\limits_{x_1}^1\dots\int\limits_{x_d}^1{\bf X}({\bf  y})\,dy_1\dots dy_d
$$
have the following representation:
$$
{\bf K}_{{\bf X}^{[0]}}={\bf T}{\bf K}_{\bf X}{\bf T}^*;\qquad {\bf K}_{{\bf X}^{[1]}}={\bf T}^*{\bf K}_{\bf X}{\bf T}.
$$
\end{lemma}

For the brevity we also introduce the notation for the integrated centered field:
$$
{\bf X}^{\{0\}}({\bf x})=\big({\overline{\bf X\vphantom{^1}}}\big)^{[0]}({\bf x}); \qquad
{\bf X}^{\{1\}}({\bf x})= \big({\overline{\bf X\vphantom{^1}}}\big)^{[1]}({\bf x}).
$$

\section{Spectral equivalence of certain Gaussian fields}

We begin with two generalizations of relation (\ref{Wa}).

Consider the stochastic integral of a non-random function $f$, see \cite[Ch.4, \S 2]{Rev},
$$
{\mathfrak F}_W(x) = \int\limits_0^x f(t) \,dW(t), \quad 0\leq x \leq 1.
$$

\begin{theorem}\label{Th1}
Let $f\in L_2(0,1)$. Then the following relation is true:
\begin{equation}
\label{FW}
\overline{{\mathfrak F}_W\vphantom{^1}}(x)\equiv {\mathfrak 
F}_W(x)-\int\limits_0^1{\mathfrak F}_W(y)\,dy \sim 
f(1-x)B(x)\sim 
f(x)B(x).
\end{equation}
\end{theorem}

\begin{proof} It is easy to check that $K_{{\mathfrak 
F}_W}=TM_f^2T^*$. So, Lemmata \ref{Le1} and \ref{Le2} 
together with relation (\ref{bridge}) imply
$$
\aligned
K_{\overline{{\mathfrak F}_W\vphantom{^1}}}= &\, \big[(I-P)TM_f\big]\cdot\big[M_fT^*(I-P)\big];\\
K_{fB}=M_fT^*(I-P)TM_f= &\, 
\big[M_fT^*(I-P)\big]\cdot\big[(I-P)TM_f\big].\\
K_{fB}=M_fT^*(I-P)TM_f= &\, [M_fT^*(I-P)TM_fS]\cdot S;\\
K_{f(1-x)B}=M_{Sf}T(I-P)T^*M_{Sf}
= &\, (SM_fS)(ST^*S)(I-P)(STS)(SM_fS)\\
= &\, S\cdot\big[M_fT^*(I-P)TM_fS\big].
\endaligned
$$
The application of Proposition \ref{Prop1} completes the 
proof.
\end{proof}

Further, we introduce a fractional counterpart of relation (\ref{Wa}). Recall that the Riemann--Liouville process is defined as the stochastic integral
$$
R_\alpha(x)=\alpha\int\limits_0^x (x-t)^{\alpha-1}\,dW(t),\qquad \alpha>\frac 12,
$$
see \cite[Sec. 3.2]{Lif} (the normalizing factor is chosen for convenience). We also define the Riemann--Liouville bridge $R_\alpha^\circ(x)$ and the generalized centered Riemann--Liouville process $R_\alpha^\bullet(x)$:
$$
R_\alpha^\circ(x)=R_\alpha(x)-x^{\alpha} R_\alpha(1);\qquad
R_\alpha^\bullet(x)=R_\alpha(x)-\alpha x^{\alpha-1} \int\limits_0^1 R_\alpha(t)\,dt.
$$

\begin{theorem}\label{Th2}
 The Riemann--Liouville bridge is spectrally equivalent to the generalized centered Riemann--Liouville process:
$$
R_\alpha^\circ(x)\sim R_\alpha^\bullet(x).
$$
\end{theorem}

\begin{proof} It is easy to check that the covariance operators of the Riemann--Liouville process has the representation $K_{R_\alpha}=T_\alpha {T_\alpha}\!\!^*$, where
$$
(T_\alpha u)(x)=\alpha\int\limits_0^x (x-t)^{\alpha-1}u(t)\,dt;\qquad
({T_\alpha}\!\!^* u)(x)=\alpha\int\limits_x^1 (t-x)^{\alpha-1}u(t)\,dt=(ST_\alpha Su)(x).
$$
We introduce the (non-orthogonal) projectors
$$
(P_\alpha u)(x)=\alpha\int\limits_0^1 (1-t)^{\alpha-1}u(t)\,dt;\qquad ({P_\alpha}\!\!^* u)(x)=\alpha(1-x)^{\alpha-1}\int\limits_0^1 u(t)\,dt
$$
and notice that $SP_\alpha=P_\alpha$, ${P_\alpha}\!\!^*S={P_\alpha}\!\!^*$.

Similarly to Lemmata \ref{Le1} and \ref{Le2}, we have
$$
\aligned
K_{R_\alpha^\circ}= T_\alpha(I-P_\alpha)(I-{P_\alpha}\!\!^*){T_\alpha}\!\!^*= &\, T_\alpha(I-P_\alpha) S^2(I-{P_\alpha}\!\!^*)ST_\alpha S\\ 
= &\,\big[T_\alpha(I-P_\alpha)S\big]\cdot\big[(I-S{P_\alpha}\!\!^*)T_\alpha S\big];
\endaligned
$$
$$
\aligned
K_{R_\alpha^\bullet}= (I-S{P_\alpha}\!\!^*)T_\alpha {T_\alpha}\!\!^*(I-P_\alpha S)= &\, (I-S{P_\alpha}\!\!^*)T_\alpha ST_\alpha S(I-P_\alpha S)\\
= &\, \big[(I-S{P_\alpha}\!\!^*)T_\alpha S\big]\cdot\big[T_\alpha(I-P_\alpha)S\big],
\endaligned
$$
and the statement follows by Proposition \ref{Prop1}.
\end{proof}

Turning to the Gaussian fields, we give a direct multivariate analog of the relation (\ref{Wa}).

\begin{theorem}\label{Th3}
 The pinned Brownian sheet is spectrally equivalent to the centered Brownian sheet:
$$
{\bf B}({\bf x})\sim \overline{\bf W\vphantom{^1}}({\bf x}).
$$
\end{theorem}

\begin{proof} For $d=1$, this relation reads $B\sim \overline{W\vphantom{^1}}$, see (\ref{Wa}). For $d=2$ it was proved in \cite{PecY} by the stochastic Fubini theorem, see (\ref{bivar}). In fact, ${\bf S}{\bf P}={\bf P}{\bf S}={\bf P}$, and Lemmata \ref{Le1} and \ref{Le2} give for any $d$
$$
\aligned
{\bf K}_{{\bf B}}={\bf T}({\bf I}-{\bf P}){\bf T}^*={\bf T}({\bf I}-{\bf P}){\bf S}{\bf T}{\bf S}=\big[{\bf T}{\bf S}({\bf I}-{\bf P})\big]\cdot\big[({\bf I}-{\bf P}){\bf T}{\bf S}\big];\\
\quad {\bf K}_{\overline{\bf W\vphantom{^1}}}=({\bf I}-{\bf P}){\bf T}{\bf T}^*({\bf I}-{\bf P})=({\bf I}-{\bf P}){\bf T}{\bf S}{\bf T}{\bf S}({\bf I}-{\bf P})=
\big[({\bf I}-{\bf P}){\bf T}{\bf S}\big]\cdot\big[{\bf T}{\bf S}({\bf I}-{\bf P})\big],
\endaligned
$$
and the statement follows by Proposition \ref{Prop1}.
\end{proof}

In the same way we derive a more general relation
\begin{equation}\label{a}
{\bf W}({\bf x})-a{\bf W}({\bf 1})\prod_{k=1}^d x_k
\sim {\bf W}({\bf x})-a\int\limits_{[0,1]^d}{\bf W}({\bf y})\,d{\bf y}, \qquad a\in\mathbb R.
\end{equation}
For $d=1$ it was proved in \cite{Dona}. If we denote the left-hand side of (\ref{a}) by ${\bf B}_a({\bf x})$ and the right-hand side by ${\bf W}_a({\bf x})$ then it is not difficult to see that
$$
{\bf K}_{{\bf B}_a}=\big[{\bf T}{\bf S}({\bf I}-a{\bf P})\big]\cdot\big[({\bf I}-a{\bf P}){\bf T}{\bf S}\big];
\quad {\bf K}_{{\bf W}_a}=\big[({\bf I}-a{\bf P}){\bf T}{\bf S}\big]\cdot\big[{\bf T}{\bf S}({\bf I}-a{\bf P})\big],
$$
and the statement follows.

\begin{rem}\label{Rem1} 
Also it is not difficult to obtain a multivariate counterpart of Theorem \ref{Th1} using the stochastic integral with respect to the Brownian sheet, see \cite{Cai}, \cite{Doz},
$$
{\mathfrak F}_{\bf W}({\bf x}) = \int\limits_0^{x_1}\dots\int\limits_0^{x_d}f({\bf  y})\,d{\bf W}({\bf y}), \quad {\bf x}\in [0,1]^d.
$$
It holds that
$$
\overline{{\mathfrak F}_{\bf W}\vphantom{^1}}({\bf x})\equiv {\mathfrak F}_{\bf W}({\bf x})-\int\limits_{[0,1]^d}{\mathfrak F}_{\bf W}({\bf y})\,d{\bf y} \sim f({\bf 1}-{\bf x}){\bf B}({\bf x}).
$$
Since ${\bf B}$ is not symmetric for $d>1$, see Remark 
\ref{Rem}, the second equivalence in (\ref{FW}) has no 
multivariate counterparts.
\end{rem}

It is appropriate here to return to the bivariate norm identities written out in \cite{DPY}. Formula (3.26) there reads:
\begin{equation}
\label{Paul}
\aligned
{\bf B}_*({\bf x}) \sim {\bf W}({\bf x}) - &\, \int\limits_0^1{\bf W}(y_1,x_2)\,dy_1 - \int\limits_0^1{\bf W}(x_1,y_2)\,dy_2\\
+ &\,  \int\limits_0^1 \int\limits_0^1 {\bf W}({\bf y})\,dy_1dy_2,\quad {\bf x}\in [0,1]^2.
\endaligned
\end{equation}

To give a new proof of (\ref{Paul}), observe that the covariance operator on the left-hand side is the tensor product of two covariance operators $K_B$. At the same time the covariance operator of the Gaussian field in the right-hand side of (\ref{Paul}) is the tensor product of two covariance operators $K_{\overline{W\vphantom{^1}}}$.  As $B$ and $\overline{W\vphantom{^1}}$ are spectrally equivalent, see (\ref{Wa}), it remains to apply our Proposition \ref{Prop2}.

Quite similarly, in the $d$-variate case we have
$$
{\bf K}_{{\bf B}_*}=
{\bf T}^*\cdot\underset{k=1}{\overset{d}\otimes}(I-P_k)\cdot{\bf T}= \big[{\bf T}^*\cdot\underset{k=1}{\overset{d}\otimes}(I-P_k)\big]\cdot\big[\underset{k=1}{\overset{d}\otimes}(I-P_k)\cdot{\bf T}\big].
$$
Using Proposition \ref{Prop1}, we obtain the field ${\bf Z}$ with the covariance operator
$$
{\bf K}_{\bf Z}=
\big[\underset{k=1}{\overset{d}\otimes}(I-P_k)\cdot{\bf T}\big]\cdot\big[{\bf T}^*\cdot\underset{k=1}{\overset{d}\otimes}(I-P_k)\big],
$$
spectrally equivalent to the $d$-variate Brownian pillow.

The direct expression of ${\bf Z}$ is more complicated and contains partial integrals with respect to all variables. For instance, the trivariate analog of (\ref{Paul}) reads
\begin{equation*}
\aligned
{\bf B}_*({\bf x}) \sim &\, {\bf W}({\bf x}) - \int\limits_0^1{\bf W}(y_1, x_2,x_3)\,dy_1 -  \int\limits_0^1{\bf W}(x_1, y_2,x_3)\,dy_2 -  \int\limits_0^1{\bf W}(x_1, x_2,y_3)\,dy_3\\
+ &\, \int\limits_0^1 \int\limits_0^1 {\bf W}(y_1, y_2,x_3)\,dy_1dy_2 + \int\limits_0^1 \int\limits_0^1 {\bf W}(y_1, x_2,y_3)\,dy_1dy_3 + \int\limits_0^1 \int\limits_0^1 {\bf W}(x_1, y_2,y_3)\,dy_2dy_3\\
- &\, \int\limits_0^1 \int\limits_0^1 \int\limits_0^1 {\bf W}({\bf y})\,dy_1dy_2dy_3,
\qquad {\bf x}\in [0,1]^3.
\endaligned
\end{equation*}

Another identity from \cite{DPY} concerns the Kiefer field ${\mathfrak K}$ whose covariance operator is the tensor product of the Wiener process and of the Brownian bridge covariances. Identity (3.27) reads
\begin{equation}
\label{Kie}
\mathfrak K({\bf x}) \sim {\bf W}({\bf x}) - \int\limits_0^1 {\bf W}(x_1,y_2)\, dy_2, \qquad {\bf x} \in [0,1]^2.
\end{equation}
To establish the operator proof of (\ref{Kie}) we observe that the covariance operator in the right-hand side is the tensor product of the covariance operators of $W$ and $\overline{W}$. It remains to apply the spectral equivalence (\ref{Wa}) and Proposition \ref{Prop2}.

\begin{rem}\label{Rem2} 
 It is not difficult to derive some weighted analogs  of obtained results similar to \cite{DPY}.
\end{rem}

\section{Integrated fields}

Now we pass to the integrated fields.

\begin{theorem}\label{Th4}
 The (left/right)-integrated pinned Brownian sheet is spectrally equivalent to the centered (left/right)-integrated Brownian sheet:
\begin{equation}
\label{Int1}
{\bf B}^{[0]}({\bf x})\sim \overline{{\bf W}^{[0]}\vphantom{^{1^1}}}({\bf x});
\qquad {\bf B}^{[1]}({\bf x})\sim \overline{{\bf W}^{[1]}\vphantom{^{1^1}}}({\bf x}).
\end{equation}
\end{theorem}

\begin{proof} We prove the first equivalence in (\ref{Int1}), the second one can be proved in the same way. For $d=1$, this relation reads $B^{[0]}\sim \overline{W^{[0]}\vphantom{^1}}$
and was discovered in \cite{BNO}. Once again, the relation ${\bf S}{\bf P}={\bf P}{\bf S}$ and Lemmata \ref{Le1} and \ref{Le2} imply that for any $d$
$$
\aligned
{\bf K}_{{\bf B}^{[0]}}={\bf T}^2({\bf I}-{\bf P}){\bf T}^{*2}={\bf T}^2({\bf I}-{\bf P}){\bf S}{\bf T}^2{\bf S}=
\big[{\bf T}^2{\bf S}({\bf I}-{\bf P})\big]\cdot\big[({\bf I}-{\bf P}){\bf T}^2{\bf S}\big];\\
{\bf K}_{\overline{{\bf W}^{[0]}\vphantom{^{1^1}}}}=
({\bf I}-{\bf P}){\bf T}^2{\bf T}^{*2}({\bf I}-{\bf P})=
\big[({\bf I}-{\bf P}){\bf T}^2{\bf S}\big]\cdot\big[{\bf T}^2{\bf S}({\bf I}-{\bf P})],
\endaligned
$$
and the statement follows by Proposition \ref{Prop1}.
\end{proof}

\begin{rem}\label{Rem3} 
 Similar relations hold for $n$-times integrated fields, for instance,
$$
{\bf B}^{[0^n]}({\bf x})\sim \overline{{\bf W}^{[0^n]}\vphantom{^{1^1}}}({\bf x}).
$$
\end{rem}

\begin{theorem}\label{Th5}
 The (left/right)-integrated pinned centered Brownian sheet is spectrally equivalent to the centered (left/right)-integrated centered Brownian sheet:
\begin{equation}
\label{Int2}
{\bf B}^{\{0\}}({\bf x})\sim \overline{{\bf W}^{\{0\}}\vphantom{^{1^1}}}({\bf x});
\qquad {\bf B}^{\{1\}}({\bf x})\sim \overline{{\bf W}^{\{1\}}\vphantom{^{1^1}}}({\bf x}).
\end{equation}
\end{theorem}

\begin{proof} We again restrict ourselves to the first equivalence in (\ref{Int2}). The relation ${\bf S}{\bf P}={\bf P}{\bf S}$, Lemmata \ref{Le1} and \ref{Le2} imply that for any $d$
$$
\aligned
{\bf K}_{{\bf B}^{\{0\}}}={\bf T}({\bf I}-{\bf P}){\bf T}({\bf I}-{\bf P}){\bf T}^*({\bf I}-{\bf P}){\bf T}^* = &\, 
\big({\bf T}({\bf I}-{\bf P})\big)^2\big(({\bf I}-{\bf P}){\bf S}{\bf T}{\bf S}\big)^2
\\
= &\, \big[\big({\bf T}({\bf I}-{\bf P})\big)^2{\bf S}\big]\cdot\big[\big(({\bf I}-{\bf P}){\bf T})^2{\bf S}\big];
\\
{\bf K\vphantom{^{1^1}}}_{\overline{{\bf W}^{\{0\}}\vphantom{^{1^1}}}}= ({\bf I}-{\bf P}){\bf T}({\bf I}-{\bf P}){\bf T}{\bf T}^*({\bf I}-{\bf P}){\bf T}^*({\bf I}-{\bf P}) = &\, 
\big(({\bf I}-{\bf P}){\bf T}\big)^2\big({\bf S}{\bf T}{\bf S}({\bf I}-{\bf P})\big)^2
\\
= &\, \big[\big(({\bf I}-{\bf P}){\bf T})^2{\bf S}\big]\cdot\big[\big({\bf T}({\bf I}-{\bf P})\big)^2{\bf S}\big],
\endaligned
$$
and the statement follows by Proposition \ref{Prop1}.
\end{proof}

\begin{rem}\label{Rem4} 
 Similar relations hold for $n$-times integrated fields, for instance,
$$
{\bf B}^{\{0^n\}}({\bf x})\sim \overline{{\bf W}^{\{0^n\}}\vphantom{^{1^1}}}({\bf x}).
$$
For $d=1$, this relation reads $B^{\{0^n\}}\sim \overline{W^{\{0^n\}}\vphantom{^1}}$ and was discovered in \cite[Sec. 4]{N09}.
\end{rem}

\begin{theorem}\label{Th6}
 The centered right-integrated pinned Brownian sheet is spectrally equi\-valent to the right-integrated centered Brownian sheet:
$$
\overline{{\bf B}^{[1]}\vphantom{^{1^1}}}({\bf x})
\sim {\bf W}^{\{1\}}({\bf x}).
$$
\end{theorem}

\begin{proof} For $d=1$ this relation reads $\overline{B^{[1]}\vphantom{^{1^1}}}\sim W^{\{1\}}$. Notice that, by symmetry of both Brownian bridge and centered Wiener process, we can change the right-integration to the left-integration that is not the case for $d>1$. The spectral equivalence $\overline{B^{[0]}\vphantom{^{1^1}}}\sim W^{\{0\}}$ was also first observed in \cite{BNO}.

The relation ${\bf S}{\bf P}={\bf P}{\bf S}$ and Lemmata \ref{Le1} and \ref{Le2} imply that for any $d$
$$
\aligned
{\bf K}_{\overline{{\bf B}^{[1]}\vphantom{^{1^1}}}}= 
({\bf I}-{\bf P}){\bf T}^*{\bf T}({\bf I}-{\bf P}){\bf T}^*{\bf T}({\bf I}-{\bf P})
= &\, \big[({\bf I}-{\bf P}){\bf S}{\bf T}{\bf S}{\bf T}({\bf I}-{\bf P}){\bf S}{\bf T}\big]\cdot\big[{\bf S}{\bf T}({\bf I}-{\bf P})\big];
\\
{\bf K\vphantom{^{1^1}}}_{{\bf W}^{\{1\}}}= {\bf T}^*({\bf I}-{\bf P}){\bf T}{\bf T}^*({\bf I}-{\bf P}){\bf T}= &\, {\bf S}{\bf T}{\bf S}({\bf I}-{\bf P}){\bf T}{\bf S}{\bf T}{\bf S}({\bf I}-{\bf P}){\bf T}\\
= &\, \big[{\bf S}{\bf T}({\bf I}-{\bf P})\big]\cdot\big[({\bf I}-{\bf P}){\bf S}{\bf T}{\bf S}{\bf T}({\bf I}-{\bf P}){\bf S}{\bf T}\big],
\endaligned
$$
and once again, the statement follows by Proposition \ref{Prop1}.
\end{proof}

\begin{rem}\label{Rem5} 
 Similar (but more intricate) relations hold for multiply integrated fields, for instance,
$$
\overline{({\bf B})^{[011]}\vphantom{^{1^1}}}({\bf x})\sim \big(\overline{{\bf W}^{[0]}\vphantom{^{1^1}}}\big)^{[11]}({\bf x}).
$$
\end{rem}

\section{Detrended processes of high order}

For the last result we restrict ourselves to the univariate case. For the Gaussian process $X$ on $[0,1]$, consider the $n$-th order detrended process (see \cite{Ai}, \cite{Pe}, and earlier papers \cite{MacN}, \cite{ALL}):
$$
X_{\langle n\rangle}(t)=X(t)-\sum\limits_{j=0}^na_jt^j,
$$
where $a_j$ are defined by the relations
$$
\int\limits_0^1 X_{\langle n\rangle}(t)t^j\,dt=0,\qquad j=0,\dots,n.
$$

\begin{theorem}\label{Th7}
 The $n$-th order detrended $n$-times integrated Wiener process is spectrally equivalent to the conditional (``bridged'') $n$-times integrated Wiener process, the Lachal process \cite{Lach}:
$$
(W^{[0^n]})_{\langle n\rangle}(t)\sim \mathbb{B}_n(t)\equiv \Big(W^{[0^n]}(t)\ \Big|\ W^{[0^m]}(1)=0,\ \ m=0,\dots,n\Big).
$$
\end{theorem}

\begin{proof} For $n=0$ this relation coincides with $\overline{W\vphantom{^1}}\sim B$. Notice that the $n$-th order detrending operation can be considered as the projection onto the subspace of $L_2([0,1])$ orthogonal to the polynomials with degree not greater than $n$. Therefore, the covariance operator $K_{(W^{[0^n]})_{\langle n\rangle}}$ has the following representation
$$
K_{(W^{[0^n]})_{\langle n\rangle}}=\big[(I-P_{\langle n\rangle})T^n\big]\cdot\big[T^{*n}(I-P_{\langle n\rangle})\big],
$$
where $P_{\langle n\rangle}$ is the orthogonal projector in $L_2([0,1])$ onto the subspace ${\cal P}_n$ of polynomials with degree not greater than $n$.

On the other hand, the direct calculation shows that the covariance operator $K_{\mathbb{B}_n}$ can be written as
$$
\aligned
K_{\mathbb{B}_n}=T^n(I-P_{\langle n\rangle})T^{*n}\stackrel{\star}{=} 
T^{*n}(I-P_{\langle n\rangle})T^n
= 
\big[T^{*n}(I-P_{\langle n\rangle})\big]\cdot\big[(I-P_{\langle n\rangle})T^n\big]
\endaligned
$$
(the relation ($\star$) holds by symmetry of $\mathbb{B}_n$), and the statement follows from Proposition \ref{Prop1}.
\end{proof}

\section{An application to goodness-of-fit tests}

Here we give an application of the obtained results to the classical goodness-of-fit problem, see, e.g., \cite{VNB}. Consider the sample $X_1,...,X_n$ with continuous distribution function $F$ in $R^d$, $d\ge2$. We are testing the simple null hypothesis $H_0$: $F=F_0$. Via the well-known Rosenblatt transform \cite{Ro} we reduce testing $H_0$ to testing uniformity on the unit cube $[0,1]^d$ using the transformed sample $\mathcal{S} = ({\bf x}^1,\dots,{\bf x}^n)$.

Let $F_n$ be the empirical distribution function based on this sample. The famous test statistic $\omega_n^2$, see \cite{SW} for its history and properties, in our case has the form
$$
\omega_n^2 = \int\limits_{[0,1]^d} \Big( F_n ({\bf z}) - \prod_{i=1}^d z_i \Big)^2 d{\bf z},
$$
which can be also written as
\begin{equation}
\label{omega}
\omega_n^2 = \Big(\,\frac 13\,\Big)^d -\frac2n \sum_{{\bf x}\in \mathcal{S}}\prod_{k=1}^d \Big(\frac{1-x_k^2}{2}\Big) + \frac{1}{n^2}\sum_{{\bf x}, {\bf x}'\in \mathcal{S}}\prod_{k=1}^d \big(1 - \max\{x_k, x'_k\}\big).
\end{equation}

 We can interpret (\ref{omega}) as a degenerate $V$-statistic which has the limiting distribution depending on the eigenvalues of its kernel
\begin{equation*}
{\cal  Q}({\bf x,y})= \prod_{k=1}^d\big(1-\max\{x_k,y_k\}\big)-2^{-d}\prod_{k=1}^d(1-x_k^2)-2^{-d}\prod_{k=1}^d(1-y_k^2)+3^{-d},
\end{equation*}
see, e.g., \cite[Ch.4]{Koro}. An old and well-known problem consists in finding the spectrum of this kernel ${\cal Q}$, i.e. the eigenvalues of the problem
\begin{equation*}\label{eigen}
({\bf Q}u)({\bf x}):=\int\limits_{[0,1]^d}{\cal Q}({\bf x,y}) u({\bf y})\,d{\bf y}=\lambda u({\bf x}),\qquad {\bf x}\in[0,1]^d,
\end{equation*}
In particular, the first eigenvalue is of special interest because it is important for the Bahadur approximate efficiency calculation of the omega-square test \cite{Bah}. It is also indispensable when evaluating the logarithmic large deviation asymptotics of $\omega_n^2$-statistic.

It is not difficult to see that
\begin{equation*}
\aligned
{\cal Q}({\bf x,y})={\cal K}_1({\bf x,y})- &\, \int\limits_{[0,1]^d}{\cal K}_1({\bf x,y})\,d{\bf y}
-\int\limits_{[0,1]^d}{\cal K}_1({\bf x,y})\,d{\bf x}\\
+ &\, \int\limits_{[0,1]^d}\int\limits_{[0,1]^d}{\cal K}_1({\bf x,y})\,d{\bf x}d{\bf y},
\endaligned
\end{equation*}
where
\begin{equation*}
{\cal K}_1({\bf x,y})=\prod_{k=1}^d\big(1-\max\{x_k,y_k\}\big)
\end{equation*}
is the covariance function of the inverted Brownian sheet ${\bf W}_1$. Therefore ${\bf  Q}$ is the covariance operator of the centered inverted Brownian sheet $\overline{{\bf W}\vphantom{^1}_1}$.

It follows from the proof of Theorem \ref{Th3} that $\overline{{\bf W}\vphantom{^1}_1}({\bf x}) \sim {\bf B}({\bf x})$. Thus, the spectrum of $\cal Q$ coincides with the spectrum of the pinned Brownian sheet ${\bf B}({\bf x})$ which was studied in many sources.

Durbin \cite{Dur} was the first to investigate the spectrum of ${\bf B}({\bf x})$ for $d=2$ and gave the list of the first $30$ reciprocal eigenvalues beginning with $\lambda_1^{-1}\approx 15.814...$. For $d=3$, the spectrum was described in \cite{Kri}, and Martynov reported the first $10$ reciprocal eigenvalues beginning by $\lambda_1^{-1}\approx 30.196...$ in \cite[\S 5]{Mart}. To the best of our knowledge, the numerical values of eigenvalues for $d>3$ are unknown.

\section*{Concluding remarks}
The referees of this paper made some suggestions which could be a line of subsequent research.

The spectral equivalence of Gaussian fields is actually much stronger than the equality in distribution of quadratic functionals. So, we can expect more applications of our results.

It would be interesting to investigate how the representations of Brownian bridges using Hardy-type operators, see \cite{JY}, \cite{GW}, \cite{Pec}, can be combined with our methods, in order to obtain new identities in distribution.

Also it is interesting to go further for various conditional $n$-times integrated Wiener process considered in \cite{Lach}.

\subsection*{Acknowledgements}
The authors are indebted to Professor I.A. Ibragimov, Professor M.A. Lifshits and anonimous referees for valuable comments and suggestions.

This work was supported by the Russian Foundation of Basic Research Grant 20-51-12004.

\end{document}